\pdfoutput=1
\documentclass[12pt]{amsart}

\usepackage[T1]{fontenc}
\usepackage{lmodern}
\usepackage{amsmath,amssymb,amsthm,microtype}
\usepackage[hidelinks]{hyperref}

\newtheorem{mainthm}{Theorem}

\newtheorem{lemma}{Lemma}[section]
\newtheorem{proposition}[lemma]{Proposition}
\newtheorem{corollary}[lemma]{Corollary}
\theoremstyle{definition}
\newtheorem{example}[lemma]{Example}
\newtheorem{remark}[lemma]{Remark}
\numberwithin{equation}{section}

\newcommand{\Pow}{\mathcal P}
\newcommand{\RR}{\mathcal R}
\newcommand{\LL}{\mathcal L}
\newcommand{\NN}{\mathbb N}
\newcommand{\ZZ}{\mathbb Z}
\newcommand{\Units}{\operatorname{U}}
\newcommand{\Aut}{\operatorname{Aut}}
\newcommand{\Sym}{\operatorname{Sym}}
\newcommand{\id}{\operatorname{id}}

\title[The Tamura--Shafer Problem for Monoids]{The Tamura--Shafer Problem for Monoids}
\author[M.~Damele]{Marco Damele}
\address{Dipartimento di Matematica e Informatica, Universit\`a degli Studi di Cagliari, Via Ospedale 72, 09124 Cagliari, Italy}
\email{marco.damele@unica.it}
\author[A.~Loi]{Andrea Loi}
\address{Dipartimento di Matematica e Informatica, Universit\`a degli Studi di Cagliari, Via Ospedale 72, 09124 Cagliari, Italy}
\email{loi@unica.it}
\date{}
\hypersetup{pdftitle={Reconstructing monoids from their power semigroups},
  pdfauthor={Marco Damele and Andrea Loi},
  pdfkeywords={power semigroups, reconstruction, global isomorphisms, monoids}}
\allowdisplaybreaks[2]
\begin{document}

\begin{abstract}
Let \(\Pow(S)\) denote the semigroup of all nonempty subsets of a semigroup \(S\), with multiplication induced from \(S\). We prove that every monoid \(M\) is reconstructible from \(\Pow(M)\), even among arbitrary semigroups: if \(\Pow(M)\cong\Pow(T)\), then \(T\) is a monoid and \(M\cong T\). The resulting isomorphism is uniquely determined by the singleton subsets of units and the two-point subsets \(\{1,a\}\) of nonunits. The main difficulty is to recover multiplication without assuming that an isomorphism of power semigroups preserves singletons. We resolve it by reconstructing finite multiplication cycles and the translation geometry of infinite mixed-power orbits. We also describe the induced split epimorphism between automorphism groups and exhibit examples where nonregular singletons are not preserved.
\end{abstract}

\maketitle

\section{Introduction}\label{sec:intro}

For a semigroup \(S\), its \emph{power semigroup} is the semigroup
\[
\Pow(S)=\{X\subseteq S:X\neq\varnothing\},\qquad XY=\{xy:x\in X,\ y\in Y\}.
\]
Although \(S\) embeds in \(\Pow(S)\) through the singleton subsets, that embedding need not be recognizable from the abstract semigroup structure of \(\Pow(S)\). The reconstruction problem is therefore not simply to restrict an isomorphism of power semigroups to the singletons: one must first find an intrinsic description of enough point labels and then recover their multiplication.

The problem originated in the work of Tamura and Shafer~\cite{TamuraShafer,Shafer}. The class of all semigroups is not determined by this construction, as Mogiljanskaja's examples show~\cite{Mogiljanskaja}. Positive results have been obtained for several substantial families, including groups~\cite{Shafer,LiTringali}, semilattices~\cite{Kobayashi}, completely \(0\)-simple semigroups~\cite{TamuraZero}, Clifford semigroups~\cite{GanZhao}, completely regular semigroups, under the generalized continuum hypothesis~\cite{YuZhao}, and commutative cancellative semigroups~\cite{Tringali}. In the inverse-semigroup setting, Liu and Tringali proved the stronger property that power-semigroup isomorphisms preserve singletons~\cite{LiuTringali}. See also~\cite{LiTringali,TringaliSurvey} for the language of global isomorphisms and further developments.

Our main result settles the reconstruction problem for all monoids, without finiteness, commutativity, cancellation, or regularity assumptions. It is deliberately formulated with an arbitrary semigroup as target; thus the identity of the source is itself reconstructed.

\begin{mainthm}\label{thm:main}
Let \(M\) be a monoid, \(T\) a semigroup, and \(\Phi:\Pow(M)\longrightarrow\Pow(T)\) an isomorphism. Then \(T\) is a monoid, and there exists a unique monoid isomorphism \(f:M\longrightarrow T\) such that
\begin{align}
\Phi(\{g\})&=\{f(g)\} &&g\in\Units(M),\label{eq:unit-label}\\
\Phi(\{1_M,a\})&=\{1_T,f(a)\} &&a\in M\setminus\Units(M).\label{eq:pair-label}
\end{align}
Moreover, \(\Phi(\{r\})=\{f(r)\}\) for every regular \(r\in M\), and \(\Phi(E)=f[E]\) for every submonoid \(E\) of \(M\). Here \(r\) is regular if \(r=rsr\) for some \(s\in M\), and \(f[E]=\{f(e):e\in E\}\).
\end{mainthm}

In particular the class of monoids is globally closed: a semigroup with the same power semigroup as a monoid is necessarily an isomorphic copy of that monoid. The last assertion of Theorem~\ref{thm:main}, applied to \(E=M\), gives \(\Phi(M)=T\). This also settles, in the monoid setting, the question about the image of the whole underlying set raised by Li and Tringali~\cite[Remark~5.5]{LiTringali}, without requiring Dedekind-finiteness of the target.

The distinction between the induced isomorphism \(f\) and the original map \(\Phi\) is essential. Even for finite monoids, \(\Phi\) need not preserve nonregular singletons. We show below that for the monoid \(M_A=\{1,0\}\sqcup A\) in which any product of two nonidentity elements is \(0\), the kernel of the natural map
\[
\Aut\Pow(M_A)\longrightarrow\Aut(M_A)
\]
is isomorphic to the full symmetric group on the nonempty subsets of \(A\) (Proposition~\ref{prop:null-kernel}). Accordingly, a proof based on universal preservation of singletons cannot work.

\subsection*{The reconstruction mechanism}
The identity of the power semigroup, when it exists, is \(\{1\}\). Its invertible elements are precisely the singletons of units. The remaining elements of \(M\) are recovered using
\[
U_a:=\{1,a\},\qquad a\notin\Units(M).
\]
The challenge is then to prove that the resulting canonical bijection is multiplicative. To explain the obstacle, observe that
\[
U_aU_b=\{1,a,b,ab\}.
\]
This is an equality of subsets, but an abstract semigroup isomorphism does not a priori preserve the inclusion relation, nor does it allow one to remove \(\{1,a,b\}\) from the product. Instead we identify \(ab\) by multiplication-theoretic tests that survive under \(\Phi\).

There are three complementary mechanisms. First, finite cyclic submonoids, infinite cyclic submonoids, powers, and regular singletons have intrinsic algebraic descriptions. Second, products of identity-containing submonoids can be inspected by equations involving \(U_x\); for suitable submonoids \(E,F\) one has
\begin{equation}\label{eq:intro-test}
x\in EF\quad\Longleftrightarrow\quad EU_xF=EF.
\end{equation}
Third, after the simpler cases are removed, one studies the mixed-power orbit
\[
O(a,b)=\{a^ib^j:i,j\ge1\}
\]
under the commuting translations \(x\mapsto ax\) and \(x\mapsto xb\). Finite cycles and noninvertible translations admit direct reconstruction. The remaining infinite orbit is a torsor for a finitely generated abelian group of rank one. Congruences arising from products \(C(a^n)C(b^m)\) recover the metric geometry of its cyclic lines, and tail tests determine the positive directions. This singles out \(ab\) even when the translation group has torsion.

\subsection*{Organization and scope}
Section~\ref{sec:labels} constructs the canonical labels, identifies cyclic submonoids, and recovers arbitrary submonoids. Section~\ref{sec:regular} treats regular elements and establishes the membership test underlying~\eqref{eq:intro-test}. Section~\ref{sec:subgroups} deals with subgroup actions, and Sections~\ref{sec:cyclic}--\ref{sec:cycles} handle unions with submonoids and finite cycles. The infinite mixed-orbit argument occupies Section~\ref{sec:mixed}; the main theorem is completed in Section~\ref{sec:completion}. Section~\ref{sec:applications} gives rigidity and automorphism-group consequences.

We emphasize that \(\Pow(M)\) here consists of \emph{all nonempty subsets}. The finitary variant and the reduced variant (finite subsets containing the identity) are distinct constructions, with their own isomorphism questions~\cite{TringaliSurvey}; counterexamples for reduced finitary power monoids are given in~\cite{Rago}. Without an identity on the source, the general finite-semigroup problem is a separate question~\cite[Problem~5.1]{LyapinNotebook}. It is still open~\cite[p.~5]{HamiltonNordahl} (see also the introduction of~\cite{LiTringali}), and Theorem~\ref{thm:main} does not address it, since the identity of the source is used throughout.

\section{Recognising cyclic submonoids and point labels}\label{sec:labels}

\subsection{Conventions and the identity}

Throughout the paper \(\NN_0=\{0,1,2,\ldots\}\). Submonoids contain the ambient identity unless they are submonoids of a local monoid with its specified identity. For a submonoid \(E\) with identity \(1\), let \(\Pow_1(E)\) denote its subsets containing \(1\). We use \(xA\) for \(\{x\}A\) and \(Ax\) for \(A\{x\}\).

The identity of \(\Pow(M)\) is \(\{1_M\}\). Conversely, if \(I\) is an identity in \(\Pow(T)\), then \(I\{t\}=\{t\}=\{t\}I\) for every \(t\in T\). Every point of \(I\) is a two-sided identity of \(T\), so \(T\) is a monoid and \(\Phi(\{1_M\})=\{1_T\}\). This observation goes back to~\cite[Lemma~1.1]{GouldIskraTsinakis}; see also~\cite{LiTringali}. In computations involving both monoids we denote either identity by \(1\) when the ambient monoid determines it.

The units of a power monoid are precisely the singletons of units. In fact, if \(XY=YX=\{1\}\), fix \(y\in Y\). For any \(x,x'\in X\),
\[
x=x(yx')=(xy)x'=x'.
\]
Interchanging \(X,Y\) also shows that \(Y\) is a singleton, and then \(xy=yx=1\), so \(x\) is a unit. Thus \(\Phi\) defines the group isomorphism \(f:\Units(M)\to\Units(T)\) of \eqref{eq:unit-label}; compare~\cite[Theorem~2.3]{LiTringali}. We extend it to nonunits below.

We will use two elementary tests for submonoids. For identity-containing submonoids \(E,F\), the relation \(E\subseteq F\) is equivalent to \(EF=F\). Also, whenever a submonoid \(D\) has only the trivial factorisation of \(1\) into two factors from \(D\), the equation \(XD=D\) recognises \(\Pow_1(D)\) in the full power monoid. Indeed, \(1\in D\) implies \(X\subseteq D\), and \(1\in XD\) then implies \(1\in X\). Conversely \(X\subseteq D\) and \(1\in X\) give \(XD=D\). We apply this test to a submonoid \(D\) of \(M\) and to its image \(\Phi(D)\): when both have only the trivial factorisation of \(1\), it gives \(\Phi(\Pow_1(D))=\Pow_1(\Phi(D))\), so that \(\Phi\) restricts to an isomorphism of reduced power monoids. In each application below this hypothesis is checked on both sides.

\subsection{Finite cyclic submonoids}

We start from a purely multiplicative way to tell whether an idempotent element of \(\Pow(M)\) is a finite submonoid containing \(1\). This is important because an isomorphism of power semigroups is not assumed to preserve cardinalities or inclusions of the underlying subsets. The counting invariant below avoids both assumptions.

\begin{lemma}\label{lem:finite-recognition}
The finite submonoids containing \(1\) are preserved by an isomorphism of power monoids. Their sizes and the numbers of ambient units in them are preserved.
\end{lemma}

\begin{proof}
For an idempotent subset \(E\) of the power monoid, define
\[
\mathcal S(E)=\{X:XE=E\},\qquad K(E)=\{x:xE\subseteq E\}.
\]
Every \(x\in K(E)\) gives a distinct element \(U_x\) of \(\mathcal S(E)\). Thus finiteness of \(\mathcal S(E)\) implies finiteness of \(K(E)\), which implies finiteness of \(E\) because \(E\subseteq K(E)\). Also \(K(E)\) is closed under products and contains \(1\). Every solution to \(XE=E\) is a subset of \(K(E)\).

Suppose \(\mathcal S(E)\) is finite. A unit in \(K(E)\) permutes \(E\) under left multiplication, since \(E\) is finite. Consequently the number \(h\) of ambient units in \(K(E)\) is the number of units of the power monoid in \(\mathcal S(E)\), and \(h\ge1\). Each subset of \(K(E)\) containing one of these units is in \(\mathcal S(E)\). Thus
\begin{equation}\label{eq:finite-count}
|\mathcal S(E)|=(2^h-1)2^{|K(E)|-h}+\ell,
\qquad 0\leq\ell<2^{|K(E)|-h},
\end{equation}
where \(\ell\) counts the solutions consisting solely of nonunits.

If \(1\in E\) then \(K(E)=E\) and \(\ell=0\). A solution must yield \(1=xy\) for some \(x\in X\), \(y\in E\). Left multiplication by \(x\) on \(E\) is surjective because \(x(yz)=z\) for \(z\in E\), hence it is injective. Comparing \(x(yx)\) with \(x1\) gives \(yx=1\), so \(X\) contains the unit \(x\). Conversely, if \(1\notin E\) then \(E\) contains no ambient unit, since finiteness and closure under powers would otherwise put \(1\) in \(E\). The solution \(X=E\) shows \(\ell>0\).

We conclude that \(E\) is a finite submonoid containing \(1\) precisely when \(E^2=E\), \(\mathcal S(E)\) is finite, and
\[
|\mathcal S(E)|=(2^h-1)2^j\quad\text{for some }j\ge0.
\]
Indeed, when \(\ell>0\), the value in \eqref{eq:finite-count} is strictly between the two successive values of this form with exponents \(|K(E)|-h\) and \(|K(E)|-h+1\). When \(\ell=0\), the test recovers \(|E|=h+j\) (here \(j\) is the \(2\)-adic valuation of \(|\mathcal S(E)|\), since \(2^h-1\) is odd), and \(h\) is the number of ambient units in \(E\). All parts of the test are preserved by \(\Phi\).
\end{proof}

\begin{lemma}\label{lem:finite-pairs}
If \(C(a)\) is finite and \(a\) is a nonunit, then \(\Phi(U_a)=U_b\) and \(\Phi(C(a))=C(b)\) for a nonunit \(b\). For a finite-order unit \(g\), \(\Phi(C(g))=C(f(g))\).
\end{lemma}

\begin{proof}
Suppose first that \(a\) is a nonunit and put \(E=C(a)\). The only ambient unit in \(E\) is \(1\). Indeed a unit among the positive powers of \(a\) would make \(a\) invertible. By Lemma~\ref{lem:finite-recognition}, the image of \(E\) is also a finite submonoid of size \(|E|\) with one unit. The family \(\Pow_1(E)\) is specified by \(XE=E\), since \(E\) is finite with trivial unit group. Indeed, in a finite monoid \(xy=1\) implies \(yx=1\), because left multiplication by \(y\) is injective and hence bijective; thus \(1\) has only the trivial factorisation in \(E\), as required by the test of Section~2.1. The corresponding statement holds for \(\Phi(E)\).

For \(X\in\Pow_1(E)\), its set powers starting at \(X^0=\{1\}\) increase by inclusion until they stabilise, and equality at one step implies equality at every subsequent step. For \(X=U_a\),
\[
U_a^i=\{1,a,\ldots,a^i\}\qquad(i\ge0).
\]
Before these sets stabilise, each step adds one point, so the first stabilisation exponent is \(|E|-1\). A set \(X\) of size at least three adds at least two points at the first step, so it stabilises at an exponent at most \(|E|-2\). The image of \(U_a\) has the same first stabilisation exponent and is distinct from \(\{1\}\). It is therefore \(U_b\) for a nonunit \(b\). Taking the stable power gives \(\Phi(E)=C(b)\).

If \(g\) is a finite-order unit, \(\Phi(C(g))\) contains \(1\) and is invariant under left multiplication by \(f(g)\). It therefore contains \(C(f(g))\), and the two sets have the same size by Lemma~\ref{lem:finite-recognition} and the group isomorphism on units.
\end{proof}

\subsection{Infinite cyclic submonoids}

Finite cyclic submonoids can be recognized by the stabilization of powers of two-point subsets. For infinite cyclic submonoids the analogous stabilization is absent. We replace it by the number of solutions of a single multiplicative equation in the reduced power monoid.

We first isolate that counting argument.

\begin{lemma}\label{lem:reduced-count}
Let \(Q\) be a monoid in which, for each \(x\ne1\), the powers \(x^i\), \(i\ge0\), are distinct. Among \(X\in\Pow_1(Q)\), the sets \(U_x\) with \(x\ne1\) are precisely the ones for which \(YX=X^3\) has exactly two solutions in \(\Pow_1(Q)\). Consequently an isomorphism between the reduced power monoids of two such monoids induces a bijection of point labels preserving powers.
\end{lemma}

\begin{proof}
When \(X=U_x\) with \(x\ne1\), a solution must contain \(1\) and \(x^2\). It lies in \(X^3\) because \(1\in X\), and it cannot contain \(x^3\) because this would give \(x^4\in YX\). Since \(x^3\in YX\) and \(x^3\notin Y\), the point \(x^3\) must be \(x^2\cdot x\), so \(x^2\in Y\). Thus its two possibilities are \(\{1,x^2\}\) and \(\{1,x,x^2\}\). If \(|X|\ge3\), take \(Y=X^2\setminus Z\) for \(Z\subseteq X\setminus\{1\}\). Every point of \(X\) belongs to \(YX\) by multiplication by \(1\) on the left, and every point of \(X^2\setminus X\) belongs to \(YX\) since it belongs to \(Y\). If \(uvw\in X^3\setminus X^2\) with \(u,v,w\in X\), then \(uv\in X^2\setminus X\), hence \(uvw\in YX\). This gives at least four distinct solutions. For \(X=\{1\}\) there is exactly one solution. This adapts to reduced power monoids the counting argument of Li and Tringali~\cite[Proposition~5.1 and Lemma~5.6]{LiTringali}, which they use for numerical monoids.

To recover powers for \(k\ge2\), the unique two-point \(Y\in\Pow_1(Q)\) with
\begin{equation}\label{eq:recover-power}
YU_x^{k-1}=U_x^{2k-1}
\end{equation}
is \(U_{x^k}\). Indeed, \(Y\subseteq C(x)\) because the right side has only powers of \(x\), and comparison of maximal exponents determines its nonidentity point.
\end{proof}

\begin{lemma}\label{lem:infinite-cyclic}
If \(C(a)\) is infinite, then \(\Phi(C(a))=C(b)\) and \(\Phi(U_a)=U_b\) for a point \(b\) with distinct powers.
\end{lemma}

\begin{proof}
Let \(C=C(a)\) and \(W=\Phi(C)\), which satisfies \(W^2=W\). Suppose \(1\notin W\), and set \(H=W\cup\{1\}\), \(B=\Phi^{-1}(H)\). Then \(BC=C=CB\) because \(HW=W=WH\). It follows that \(B\subseteq C\). As \(1\in BC\) and \(C\) is a free monogenic monoid, \(1\in B\). Thus \(B\) and \(H\) are both submonoids, each with only the trivial factorisation of \(1\) within itself. For \(H\), this follows from the closure of \(W\) under products. The tests \(XB=B\) and \(XH=H\) identify their reduced power monoids under \(\Phi\).

Any nonidentity element of \(H\) with finitely many powers would generate a nontrivial finite submonoid in \(H\). By Lemma~\ref{lem:finite-recognition} and the inclusion test for submonoids, its preimage would be a nontrivial finite submonoid in \(B\subseteq C\). Thus the powers of each nonidentity element of \(H\) are distinct. The same argument applies to any identity-containing submonoid matched under \(\Phi\) to a submonoid of \(C\).

Apply Lemma~\ref{lem:reduced-count} to \(B,H\). Any two nonidentity elements of \(B\) have a common positive power, so the same is true in \(H\). In \(H\), if \(x,y\neq1\) and \(xy=x\), then \(xy^s=x\) for all \(s\ge1\). Choosing \(y^s=x^r\) with \(r,s>0\) gives \(x^{r+1}=x\), which is impossible. The dual argument excludes \(xy=y\).

Let \(b_0\) be the element of \(B\setminus\{1\}\) with least exponent as a positive power of \(a\). Such a point exists because \(H\ne\{1\}\). Removing \(b_0\) leaves a submonoid, since a product of two nonidentity points of \(B\) has larger exponent. Within each reduced power monoid, submonoids are the idempotents and membership \(x\in E\) is recognised by \(U_xE=E\). Thus removing the corresponding point \(h_0\) from \(H\) also leaves a submonoid. Indeed, the isomorphism of reduced power monoids maps the submonoid \(B\setminus\{b_0\}\) to a submonoid of \(H\), and by the membership test its points are exactly the labels corresponding to the points of \(B\setminus\{b_0\}\), that is, all points of \(H\) except \(h_0\). But \(h_0\in W=W^2\), so \(h_0=xy\) with \(x,y\in W\). Both factors differ from \(h_0\) by the absence of absorption. This contradicts the closure of \(H\setminus\{h_0\}\).

Hence \(1\in W\) and the powers of each nonidentity element of \(W\) are distinct. If \(w\) is a unit of \(W\), the ambient unit singleton \(\{w\}\) satisfies \(\{w\}W=W\). Its preimage is an ambient unit singleton \(\{g\}\) satisfying \(\{g\}C=C\). This implies both \(g\in C\) and that its inverse belongs to \(C\), giving \(g=1\). If \(xy=1\) in \(W\) then \(yx\) is idempotent, so \(yx=1\) by distinctness of powers for nonidentity elements. Thus \(x,y\) are units of \(W\), giving \(x=y=1\). It follows that \(XW=W\) recognises \(\Pow_1(W)\). Lemma~\ref{lem:reduced-count} now matches the point labels of \(C,W\) and all their powers. In particular the point matching \(a\) is a generator \(b\) of \(W\) as required.
\end{proof}

\subsection{Aligning unit and nonunit labels}

A priori the preceding constructions give two systems of labels: singleton labels for invertible elements and two-point labels for the remaining elements. An infinite-order unit also belongs to an infinite cyclic submonoid, so its two possible labels must be shown to agree. This is the purpose of the next lemma.

\begin{lemma}\label{lem:submonoids}
Equations~\eqref{eq:unit-label} and \eqref{eq:pair-label} define a bijection \(f:M\to T\). Moreover,
\[
\Phi(C(x))=C(f(x)),\qquad \Phi(E)=f[E]
\]
for \(x\in M\) and every submonoid \(E\) of \(M\).
\end{lemma}

\begin{proof}
First let \(g\) be an infinite-order unit and put \(g'=f(g)\), \(H=\langle g\rangle_{\rm gr}\), and \(H'=\langle g'\rangle_{\rm gr}\). Here the notation \(\langle\cdot\rangle_{\rm gr}\) denotes generation as a group. Since \(H=C(g)C(g^{-1})\), the image \(\Phi(H)\) contains \(1\) by Lemma~\ref{lem:infinite-cyclic}. The equalities \(\{g\}H=H=\{g^{-1}\}H\) give \(H'\subseteq\Phi(H)\). In the other direction the same reasoning gives \(H\subseteq\Phi^{-1}(H')\). Both inclusions transport by the submonoid inclusion test, so \(\Phi(H)=H'\).

Each of the rays \(C(g),C(g^{-1})\) is maximal among proper submonoids of \(H\), because adjoining any point of the other ray yields a generator of the other sign. Their images are infinite cyclic submonoids of \(H'\). If one were generated by a proper power of \(g'\) or \((g')^{-1}\), it would be properly contained in one of the two generator rays. The preimage of that ray contains \(1\) by Lemma~\ref{lem:infinite-cyclic}, contradicting maximality. The image rays have opposite signs since their product is \(H'\). Write \(\Phi(U_g)=U_b\) and \(R=\Phi(C(g^{-1}))\). Then \(b\) generates the opposite ray to \(R\). In the source,
\[
U_g C(g^{-1})=\{g\}C(g^{-1}).
\]
In the target \(U_bR=\{b\}R\), so \(\{b\}R=\{g'\}R\). Distinct translates of a ray in an infinite cyclic group are distinct. It follows that \(b=g'\).

By Lemmas~\ref{lem:finite-pairs} and \ref{lem:infinite-cyclic}, the image of a nonunit pair is a pair, and its nonidentity point is a nonunit. In the infinite cyclic case an infinite-order unit image is ruled out by applying the preceding alignment to \(\Phi^{-1}\). The same arguments for \(\Phi^{-1}\) give a bijection on the nonunit labels. Together with the unit labels this gives \(f\) and the correspondence for cyclic submonoids.

For any two idempotents \(F,E\) in a semigroup, the natural partial order is \(F\preceq E\) if \(FE=EF=F\). Let \(E\) be the submonoid generated by a family of elements \(x\). Within the idempotents of the power semigroup, \(E\) is the greatest lower bound of their cyclic submonoids \(C(x)\). Indeed \(E\preceq C(x)\) for every generator. If an idempotent subset \(F\) satisfies \(F\preceq C(x)\) for each generator, then \(Fx\subseteq F\) and \(xF\subseteq F\). Iterating on words in the generators gives \(FE=EF=F\), using \(1\in E\). Thus \(F\preceq E\) even when \(F\) does not contain \(1\).

Thus \(\Phi(E)\) is the submonoid generated by the \(C(f(x))\) for \(x\in E\). For a nonunit \(x\), membership in an identity-containing submonoid is recognised by \(U_x E=E\). For a unit it is recognised by \(C(x)\subseteq E\). These tests give \(\Phi(E)=f[E]\).
\end{proof}

\begin{remark}\label{rem:transport}
We say that a subset \(X\subseteq M\) is \emph{recognised}, or \emph{matched pointwise}, if \(\Phi(X)=f[X]\). By the definition of \(f\), the singletons of units and the pairs \(U_a\) of nonunits are recognised, and Lemma~\ref{lem:submonoids} adds all submonoids. Since \(\Phi\) is an isomorphism, an equation between products of recognised subsets holds in \(\Pow(M)\) if and only if the corresponding equation between their images holds in \(\Pow(T)\); this is what we mean when we say that a test \emph{transports}. The list of recognised subsets grows along the paper: regular singletons (Proposition~\ref{prop:regular}), the members of \(\Pow_1(C(x))\) for infinite \(C(x)\) (Lemma~\ref{lem:powers}), products \(EF\) of submonoids of \(D_R\) or \(D_L\) (Lemma~\ref{lem:two-products}), and the subgroups and corners used in Section~\ref{sec:subgroups}.
\end{remark}

\begin{example}[Units of finite order]\label{ex:units}
For units of finite order, the pair labels cannot replace the singleton labels in~\eqref{eq:unit-label}. Let \(G=\ZZ/3=\{1,g,g^2\}\). For a fixed \(k\in G\), the bijection of \(\Pow(G)\) that fixes the three singletons and \(G\), and sends each co-singleton \(G\setminus\{x\}\) to \(G\setminus\{xk\}\), is an automorphism: it respects the products \(\{h\}(G\setminus\{x\})=G\setminus\{hx\}\) because \(G\) is abelian, and the product of two co-singletons is \(G\). For \(k=g\) it maps \(U_g=G\setminus\{g^2\}\) to \(G\setminus\{1\}=\{g,g^2\}\), which does not even contain \(1\). The same example shows that the parameter \(h_0\) in the proof of Proposition~\ref{prop:group-action} cannot be avoided, and that the kernel \(K(G)\) of Corollary~\ref{cor:aut} can be nontrivial for a group, although all singletons are preserved by Corollary~\ref{cor:regular-generated}. Here \(|\Aut(\Pow(G))|=6\) and \(K(G)\cong\ZZ/3\).
\end{example}

\section{Idempotents, regular elements and products of submonoids}\label{sec:regular}

\begin{lemma}\label{lem:fixed-sets}
For a point \(x\), the nonempty subsets fixed pointwise by multiplication on the left by \(x\) are precisely the subsets \(X\) such that \(Y=X\) is the unique solution to \(U_xY=X\). The analogous statement holds on the right.
\end{lemma}

\begin{proof}
If \(X\) is fixed pointwise on the left, any solution \(Y\subseteq X\) satisfies \(U_xY=Y\), so \(Y=X\). Conversely, if \(Y=X\) is a solution, then \(xX\subseteq X\). If \(t\in X\) has \(xt\ne t\), the nonempty set \(Y=X\setminus\{xt\}\) still gives \(U_xY=X\), because \(t\) remains in \(Y\) and supplies \(xt\). This contradicts uniqueness. The right-handed proof is the same with reversed products.
\end{proof}

\begin{proposition}\label{prop:regular}
The isomorphism \(\Phi\) sends each regular singleton \(\{r\}\) to \(\{f(r)\}\). The map \(f\) preserves products with an idempotent factor, and it preserves a product with a regular factor whenever the product is regular.
\end{proposition}

\begin{proof}
If \(e\ne1\) is idempotent, \(U_e^2=U_e\) and preservation of the nonunit pair show that \(f(e)\) is idempotent. Applying both tests of Lemma~\ref{lem:fixed-sets} to \(e\) identifies the nonempty subsets of \(eMe\) and, in the target, those of \(f(e)Tf(e)\). Their identity elements are the singletons of \(e\) and \(f(e)\), respectively. Thus \(\Phi(\{e\})=\{f(e)\}\), including \(e=1\).

For a proper idempotent \(e\), the submonoid \(\{1\}\cup eM\) is the least submonoid \(F\) satisfying
\[
U_e F=F,\qquad \{e\}F=\{e\}M.
\]
Indeed the conditions imply \(eM=eF\subseteq F\), and \(\{1\}\cup eM\) satisfies both. Since \(M\) is matched to \(T\) as a subset by Lemma~\ref{lem:submonoids}, we recognise \(eM\) after adjoining \(1\). If \(x\) is a nonunit, then \(ex\) is a nonunit for \(e\ne1\) and \(f(ex)\) lies in \(f(e)T\). The equality \(\{e\}U_x=\{e\}U_{ex}\) therefore gives
\[
\{f(e),f(e)f(x)\}=\{f(e),f(ex)\}.
\]
This gives \(f(ex)=f(e)f(x)\), also when both products equal \(f(e)\). Reversing products treats \(xe\).

Let \(r\) be regular. Choose \(s\) with \(rsr=r\) and \(srs=s\), called an inverse of \(r\). These two conditions can be arranged by replacing any witness \(s_0\) to \(rs_0r=r\) with \(s_0rs_0\). Both \(rs\) and \(sr\) are idempotent. For \(A=\Phi(\{r\})\) and \(B=\Phi(\{s\})\) we have \(AB=\{p\}\), \(BA=\{q\}\) with \(p=f(rs)\), \(q=f(sr)\). Also \(pA=A=Aq\). Since \(p\) and \(q\) are idempotent, these equalities hold pointwise: if \(x=x'q\) with \(x'\in A\), then \(xq=x'q^2=x\), and similarly \(px=x\). For \(x,x'\in A\) and \(y\in B\) this gives
\[
x=xq=x(yx')=(xy)x'=px'=x'.
\]
Hence \(\Phi(\{r\})=\{h(r)\}\) with \(h(r)\) regular. The labels \(h(r), f(r)\) already agree for units and idempotents. Otherwise \(s\) is a nonunit, and, putting \(e=rs\), we have \(\{r\}U_s=\{e\}U_r\). Applying \(\Phi\) gives \(h(r)=f(e)f(r)=f(r)\), since \(h(r)\neq f(e)\) and the only other point on the right is \(f(e)f(r)\). The previous idempotent calculation applies here because \(r\) is a nonunit.

Products of two regular elements with regular output now follow from their singleton products. In particular \(eg,ge\) for a unit \(g\) are regular, so this proves the idempotent assertion for the remaining factors. For the last assertion, suppose \(r\) is a regular left factor with inverse \(s\), and \(y=rx\) is regular. Put \(e=sr\). Then \(ex=sy\) and \(y=r(ex)\), so \(ex\) and \(y\) generate the same principal left ideal. Regularity is constant on this class. Indeed if \(z=uy\), \(y=vz\) and \(y=ywy\), then \(z=z(wv)z\). Since \(ex,r,y\) are all regular, we get
\[
f(y)=f(r)f(ex)=f(r)f(e)f(x)=f(r)f(x).
\]
The same argument in the opposite monoid treats a regular right factor.
\end{proof}

\begin{lemma}\label{lem:powers}
For all \(x\in M\) and \(k\ge0\), \(f(x^k)=f(x)^k\). For infinite \(C(x)\), the restriction of \(\Phi\) to \(\Pow_1(C(x))\) is pointwise. For finite \(C(x)\) with \(x\) a nonunit, membership of any power before the periodic part in each \(X\in\Pow_1(C(x))\) is preserved.
\end{lemma}

\begin{proof}
The power assertion for units follows from the group isomorphism. For infinite \(C(x)\), use Lemmas~\ref{lem:infinite-cyclic} and \ref{lem:reduced-count}. They give \(\Phi(U_{x^k})=U_{f(x)^k}\), that is, \(f(x^k)=f(x)^k\); the pointwise assertion on \(\Pow_1(C(x))\) follows from the membership test~\eqref{eq:tail-membership} at the end of the proof. Suppose \(A=C(x)\) is finite with \(x\) a nonunit. We first show how to recognise its right ideals. Given a submonoid \(E\subseteq A\) and \(I=E\setminus\{1\}\), the tests are
\begin{equation}\label{eq:finite-right-ideals}
\bigl|\{EX:X\in\Pow_1(A)\}\bigr|=2^{|A\setminus E|},
\qquad eA\subseteq E\quad(e\in I,\ e^2=e).
\end{equation}
If \(I\) is a right ideal, \(EX=E\cup X\), so both hold. Conversely every \(EX\) lies between \(E\) and \(A\). The count forces every set \(E\cup\{t\}\), \(t\in A\setminus E\), to occur. Its preimage \(X\) contains \(t\), since \(X\subseteq EX\) and \(EX=E\) for \(X\subseteq E\). Thus \(It\subseteq I\cup\{t\}\). If \(it=t\) for \(i\in I\), choose a positive idempotent power \(e\) of \(i\). Then \(e\in I\) and \(et=t\), contradicting the second test. This gives \(It\subseteq I\), which suffices since \(IE\subseteq I\). The tests transport by recognition of \(\Pow_1(A)\), Lemma~\ref{lem:submonoids}, and preservation of products with idempotents.

The principal right ideals in \(A\) are now recovered as the least right ideals containing their generators. We have \(x^iA=\{x^j:j\ge i\}\), a strictly descending list until the periodic part begins. That part is the cyclic group of eventual powers, with identity the unique positive idempotent power \(e\). To see this, let \(r\) be the period of the repeating powers. A sufficiently large multiple of \(r\) is an exponent of \(e\), and exponents in the repeating range whose sum is a multiple of \(r\) give inverse pairs. The periodic part is \(eA\), so it is matched under \(f\) by Proposition~\ref{prop:regular}. The strict order of the right ideals recovers the earlier powers one by one and determines the first periodic index. Explicitly, let \(m\) be the first periodic index. The right ideals of \(A\) not containing \(1\) are the sets \(x^iA\), and they form the chain \(xA\supsetneq x^2A\supsetneq\cdots\supsetneq x^mA\). Since the tests above transport, \(f\) maps this chain bijectively and inclusion-preservingly onto the corresponding chain for \(C(f(x))\). Hence \(f[x^iA]=f(x)^iC(f(x))\) for \(1\le i\le m\), and since \(x^iA\setminus x^{i+1}A=\{x^i\}\) for \(i<m\), we get \(f(x^i)=f(x)^i\) for \(i<m\). The cyclic group multiplication and \(f(ex)=f(e)f(x)\) give, for \(k\) in the periodic range,
\[
f(x^k)=f((ex)^k)=f(ex)^k=(f(e)f(x))^k=f(x)^k.
\]

To recover membership, put \(E_n(x)=\{1\}\cup\{x^i:i\ge n\}\) for \(n\ge1\). For \(X\in\Pow_1(C(x))\), \(E_n(x)X=E_n(x)\cup X\). Whenever \(x^n\) precedes the periodic part, or whenever \(C(x)\) is infinite, \(E_n(x)\setminus E_{n+1}(x)=\{x^n\}\). Therefore
\begin{equation}\label{eq:tail-membership}
x^n\in X\quad\Longleftrightarrow\quad E_n(x)X=E_{n+1}(x)X.
\end{equation}
This test uses only recognised submonoids and proves the assertions about membership. For the rigidity of the reduced power monoid of \(\NN_0\) as an abstract monoid, see Tringali and Wen~\cite{TringaliWen}.
\end{proof}

\begin{lemma}\label{lem:two-products}
Define
\[
\begin{aligned}
D_R&=\{1\}\cup\{x:x\text{ has no right inverse in }M\},\\
D_L&=\{1\}\cup\{x:x\text{ has no left inverse in }M\}.
\end{aligned}
\]
Each is a submonoid with only the trivial factorisation of \(1\), and \(f\) matches them to their counterparts in \(T\). If \(D\) is either one and \(E,F\subseteq D\) are submonoids, then \(\Phi(EF)=f[EF]\).
\end{lemma}

\begin{proof}
An element with no right inverse retains that property after multiplication on the right. The dual property holds for elements with no left inverse. Hence both sets are submonoids. If \(xy=1\) with \(x,y\in D_R\), then \(x=1\) and \(y=1\). In \(D_L\) use the left inverse of \(y\) instead. Moreover any two witnesses to \(xy=1\) are regular, since \(xyx=x\) and \(yxy=y\). Thus the tests for one-sided inverses transport by Proposition~\ref{prop:regular}.

The image \(\Phi(EF)\) is the product of the corresponding submonoids, so it lies in the corresponding \(D\). For \(z\in D\),
\begin{equation}\label{eq:submonoid-product}
z\in EF\quad\Longleftrightarrow\quad E U_z F=EF .
\end{equation}
If \(z=ef\) with \(e\in E\), \(f\in F\), then \(EzF\subseteq EF\) by closure of \(E,F\) under products. Conversely \(z\in E U_z F\) because \(1\in E\cap F\). Every point \(z\in D\) has its pair recognised, including \(U_1=\{1\}\). This is where the restriction to a common \(D\) is used: every point \(z\neq1\) of \(D\) is a nonunit, whereas the pair of a unit of finite order need not be recognised (Example~\ref{ex:units}). Thus \eqref{eq:submonoid-product} tests membership in \(EF\) and its image.
\end{proof}

\section{Multiplication by subgroup elements}\label{sec:subgroups}

An idempotent \(e\) is the identity of the local monoid \(eMe\). Any subgroup of \(M\) with identity \(e\) is a subgroup of \(\Units(eMe)\), and all its elements are regular. Proposition~\ref{prop:regular} handles products of regular elements with regular output. The remaining task in this section is to determine the action of a subgroup on nonregular elements. We work with cyclic subgroups because their orbits and stabilisers suffice to test a single multiplier.

\begin{lemma}\label{lem:co-singleton}
Let \(G\) be a group with at least three elements. A nonempty \(V\subseteq G\) is the complement of a singleton precisely when, for nonempty \(Y\subseteq G\),
\[
VY=G\quad\Longleftrightarrow\quad |Y|\ge2.
\]
In that case \(V^2=G\).
\end{lemma}

\begin{proof}
For the complement of a singleton, distinct right translates omit distinct points, and the union of two such translates is \(G\). For a proper \(V\) that omits distinct \(s,t\), the two-set \(Y=\{1,t^{-1}s\}\) gives \(s\notin VY\). Also \(V=G\) fails the stated condition on singleton \(Y\). The assertion about \(V^2\) follows since \(|V|\ge2\).
\end{proof}

\begin{proposition}\label{prop:group-action}
For any \(g\in M\) belonging to a subgroup of \(M\) and any \(a\in M\),
\[
f(ag)=f(a)f(g),\qquad f(ga)=f(g)f(a).
\]
\end{proposition}

\begin{proof}
We prove the right-handed formula. Let \(H=\langle g\rangle_{\rm gr}\) have identity \(e\). Replacing \(a\) by \(ae\) is legitimate because \(ag=(ae)g\), and \(f(ae)=f(a)f(e)\) by Proposition~\ref{prop:regular}. Thus assume \(ae=a\). If \(a\) is regular, then \(ag\) generates the same principal right ideal as \(a\), using \((ag)g^{-1}=a\). Regularity is constant on that class by the right-handed version of the calculation in Proposition~\ref{prop:regular}. The formula follows from that proposition. Henceforth \(a\) is nonregular.

\medskip\noindent
\emph{The orbit and its stabiliser.}
Put \(O=aH\) and \(K=\{h\in H:ah=a\}\). The set \(H\) is preserved pointwise by the formula \(H=\{e\}(\{1\}\cup H)\{e\}\), and so is each subgroup of \(H\). The equation \(XH=H\) in the corner \(\Pow(eMe)\) recognises \(\Pow(H)\), since \(H\) contains \(e\) and is a group. Among nonregular \(x\) with \(xe=x\), the equation
\begin{equation}\label{eq:orbit-recognition}
U_xH=U_aH
\end{equation}
recognises \(O\). Indeed, \(xH\) consists of nonregular elements, because all its points generate the same principal right ideal as \(x\). Thus \(U_xH=H\cup xH\), and \(x\in xH\). The image of \(O\) as a set of labels is \(O'=f(a)H'\), where \(H'=f[H]\). To identify \(K\), for each \(h\in H\) use the same test with the cyclic subgroup generated by \(h\). The orbit has size one exactly when \(ah=a\) or, in the target, \(f(a)f(h)=f(a)\). This gives \(K'=f[K]\).

The group \(H\) is cyclic, so \(K\) is normal. The sets in \(\Pow(H)\) satisfying \(XK=X\) are exactly the nonempty unions of cosets of \(K\). Their multiplication is the power monoid of \(H/K\) under the natural quotient map. Consequently \(\Phi\) identifies the power monoids of \(H/K\) and \(H'/K'\) (here we use \(\Phi(K)=f[K]=K'\), as the subgroups of \(H\) are preserved pointwise), and the co-singleton test of Lemma~\ref{lem:co-singleton} applies there. The distinction between singleton and larger subsets within a group power monoid is provided by its units.

\medskip\noindent
\emph{Fixing the part of a set in the group.}
Let \(L_e\) be the set of elements of \(eMe\) without a left inverse in \(eMe\), and put \(Q_e=\{e\}\cup L_e\). The defining one-sided inverse equations have regular witnesses and are preserved by Proposition~\ref{prop:regular}. Thus \(Q_e\) is preserved pointwise as a subset by applying Lemma~\ref{lem:submonoids} to \(\{1\}\cup Q_e\) and multiplying by \(\{e\}\) on both sides.

For \(y\in O\), \(ey\in L_e\). Indeed \(k(ey)=e\) with \(k\in eMe\) would give \(ky=e\) and hence \(yky=ye=y\). Multiplication by \(H\) on the right preserves \(L_e\), because all its multipliers are units of \(eMe\). For \(X\subseteq H\cup O\) with \(X\cap H\neq\varnothing\), it follows that
\begin{equation}\label{eq:unit-part}
Q_e\{e\}X=(X\cap H)\cup L_e.
\end{equation}
Here the part from \(X\cap H\) gives all of \(L_e\) as well as the indicated points of \(H\). The part from \(O\) is contained in \(L_e\) since \(L_e\) is a left ideal of \(eMe\). If \(X\cap H=\varnothing\) the product is contained entirely in \(L_e\). The calculation applies to any further points \(y\) with \(ey\in L_e\) added to the allowed support of \(X\).

\medskip\noindent
\emph{A finite orbit.}
When \(|O|=1\) the generator fixes the only point. When \(|O|=2\) the generator interchanges the two points. Orbit recognition gives the corresponding action in both cases. Let \(|O|\ge3\) and set
\[
N=\bigcup_{n\ge2}O^n,\qquad E=\{1\}\cup O\cup N,\qquad E_N=\{1\}\cup N.
\]
We need \(O\cap N=\varnothing\). Suppose \(o_1\cdots o_n\in O\), \(n\ge2\). Write \(o_n=ah_n\) and \(o_1\cdots o_n=ah_*\) with \(h_n,h_*\in H\). Then \(w a=ah_*h_n^{-1}\) for \(w=o_1\cdots o_{n-1}\in aM\), since \(ae=a\). If \(h=h_*h_n^{-1}\), iteration gives \(w^j a=ah^j\). Since \(O\) is finite, some positive power of \(h\) fixes \(a\). Then \(a=w^j a\) with \(w^j\in aM\), contradicting nonregularity.

If \(e\ne1\) then \(1\notin O\cup N\) because all points of \(O\cup N\) are fixed on the right by \(e\). If \(e=1\), a product of points of \(O\) cannot be \(1\) because its last factor has no left inverse. Also \(NH=N\), since the last factor runs through \(O\). Finally
\[
e(o_1\cdots o_n)=(eo_1)\cdots(eo_n)\in L_e
\]
uses \(o_i e=o_i\) and the fact that \(L_e\) is a left ideal. Thus \(eN\subseteq L_e\) and \(N\cap H=\varnothing\).

The submonoid \(E\) is generated by \(O\), so its image is the submonoid generated by \(O'\). The disjointness calculation in the two monoids shows that the set of labels \(N'=f[N]\) is exactly \(\bigcup_{n\ge2}(O')^n\). Let \(V=H\setminus K\). The complement test in the quotient gives \(\Phi(V)=V'=H'\setminus K'h_0\) for some \(h_0\in H'\). Since the quotient has size \(|O|\ge3\), \(V^2=H\) and \((V')^2=H'\).

The sets of the form \(V\cup N\cup B\), \(B\subseteq O\), are specified by
\begin{equation}\label{eq:finite-orbit-families}
X\{e\}=X,\qquad EX=EV,\qquad Q_e\{e\}X=Q_eV,\qquad E_NX=X .
\end{equation}
The second equation confines \(X\) to \(EV=V\cup N\cup O\) (so the first equation is in fact implied by the second, since every point of \(EV\) is fixed on the right by \(e\)). The third uses \eqref{eq:unit-part} to require its part in \(H\) to be \(V\). The fourth then requires \(N\subseteq X\) since \(NV=N\). Conversely every such set satisfies the four equations, using that products of \(N\) with \(N\) or \(O\) on either side lie in \(N\).

For \(y\in O\) write \(X_y=V\cup N\cup\{y\}\). We have
\begin{equation}\label{eq:finite-orbit-two-tests}
U_y X_y=EV,\qquad X_y V=H\cup N\cup(O\setminus\{y\}).
\end{equation}
Here \(yV\) omits just \(y\) from \(O\), and \(yN\) and \(y^2\) contribute only to \(N\). In the target write \(\Phi(X_y)=V'\cup N'\cup B'\) by \eqref{eq:finite-orbit-families}. The first test in \eqref{eq:finite-orbit-two-tests} forces \(f(y)h_0\in B'\), since \(f(y)V'\) omits precisely \(f(y)h_0\). If \(B'\) contained two distinct points then \(B'V'=O'\), contradicting \(X_yV\ne EH\). Hence \(B'=\{f(y)h_0\}\). The only point of \(O'\) omitted by \(\Phi(X_yV)\) is therefore \(f(y)h_0^2\). In the source \(X_{yg}V=(X_yV)\{g\}\) because right multiplication by \(g\) fixes \(H,N\) setwise and moves the omitted point from \(y\) to \(yg\). In the target this equality gives
\[
f(yg)h_0^2=f(y)h_0^2 f(g).
\]
Commutativity in the cyclic group \(H'\) and cancellation by its action give \(f(yg)=f(y)f(g)\).

\medskip\noindent
\emph{An infinite orbit.}
Here \(H\) is infinite cyclic and \(K=\{e\}\). Put \(V=H\setminus\{e\}\), so \(\Phi(V)\) is a co-singleton in \(H'\). Fix a nonempty \(W\subseteq H\). The sets \(W\cup B\) with \(\varnothing\ne B\subseteq O\) are recognised by
\[
X\{e\}=X,\qquad XH=U_aH,\qquad Q_e\{e\}X=Q_eW.
\]
Indeed the first two equations say \(X\subseteq H\cup O\) and force \(X\) to meet both parts. The third gives its precise part in \(H\) by \eqref{eq:unit-part}. For \(|W|\ge2\), \(WV=H\) by the co-singleton test. The condition \(XV\ne H\cup O\) then detects that \(B\) has one point, since any two distinct translates on a free cyclic orbit of the co-singleton miss different points. The same holds in the target. In particular \(\Phi\) induces a bijection \(h_W:O\to O'\) from the subsets \(W\cup\{y\}\) to those of the form \(\Phi(W)\cup\{h_W(y)\}\), using \(|\Phi(W)|\ge2\) in \(\Pow(H')\).

For such \(W\), the equality \((W\cup\{y\})V=(H\cup\{y\})V\) shows that \(h_W=h_H\), since the omitted point from the orbit distinguishes the translates of a co-singleton. The map \(h_H\) intertwines right multiplication by \(g\) and by \(f(g)\), as is seen by multiplying \(H\cup\{y\}\) by \(\{g\}\). Take now \(W=\{e,g,g^2,\ldots\}\). This ray is matched pointwise as a subset by applying Lemma~\ref{lem:submonoids} to \(\{1\}\cup W\) and restricting to the corner. The equality
\[
(W\cup\{y\})W=U_y W
\]
gives \(h_H(y)\Phi(W)=f(y)\Phi(W)\) on the \(O'\)-part. The action is free and distinct translates of a ray on an infinite cyclic orbit are distinct. Thus \(h_H=f\) on \(O\), proving the assertion. The left-handed formula follows by reversing both multiplications.
\end{proof}

\section{Unions with submonoids and cyclic products}\label{sec:cyclic}

From now on we transport the target law to the set \(M\) using the point labels \(f\). Thus \(x\star y=f^{-1}(f(x)f(y))\), and the power-monoid isomorphism is \(\Psi(X)=f^{-1}[\Phi(X)]\) from the original law on subsets to the \(\star\)-law. All units as singleton subsets, nonunit pairs, individual powers, and submonoids have the same point labels. An equality involving arbitrary subsets still transports by applying \(\Psi\) to the factors. We write \(D\) for either \(D_R\) or \(D_L\), which have the same point labels in the two laws by Lemma~\ref{lem:two-products}. The test \(XD=D\) recognises \(\Pow_1(D)\), and the corresponding test works for each submonoid of \(D\).

\begin{lemma}\label{lem:union}
For submonoids \(E,S\subseteq D\), the condition
\[
EX=E\cup X\qquad\text{for every }X\in\Pow_1(S)
\]
holds in one law if and only if it holds in the other. The analogous statement holds for \(XE=E\cup X\).
\end{lemma}

\begin{proof}
It suffices to establish \(E\star U_t=E\cup U_t\) for each \(t\in S\), since multiplication by a set takes the union of the products with its individual points. Fix \(t\) and put \(A=C(t)\). The condition \(EA=E\cup A\) is transported pointwise by Lemma~\ref{lem:two-products}. If \(t\in E\), closure of \(E\) gives the claim, so assume \(t\notin E\).

When \(A\) is finite, the source condition gives
\[
\bigl|\{EX:X\in\Pow_1(A)\}\bigr|=2^{|A\setminus E|}.
\]
In the other law, each output is a set between \(E\) and \(E\cup A\). Since \(\Pow_1(A)\) is recognised, the count gives in particular an \(X\) with \(E\star X=E\cup\{t\}\). The input \(X\) contains \(t\), because \(X\subseteq E\star X\) and \(E\star X=E\) when \(X\subseteq E\). It follows that \(E\star t\subseteq E\cup\{t\}\).

When \(A\) is infinite, suppose \(e\star t=t^k\notin E\cup\{t\}\) for \(e\in E\). This is the only possible form of a new point because \(E\star A=E\cup A\), and it has \(k\ge2\). Closure of \(E\) gives \(E\star t^k\subseteq E\star t\), so \(E\star\{1,t,t^k\}=E\star U_t\). Both inputs are matched pointwise by Lemma~\ref{lem:powers}. But the equality fails in the original law, since multiplication by \(E\) adds \(E\) by union and \(t^k\notin E\cup\{t\}\). Thus again \(E\star t\subseteq E\cup\{t\}\), giving the desired formula. Interchanging the laws proves equivalence and reversing multiplication proves the right-sided version.
\end{proof}

\begin{lemma}\label{lem:absorption}
If \(a,b\) lie in a common \(D_R\) or \(D_L\) and \(ab\in\{a,b\}\), then \(a\star b=ab\).
\end{lemma}

\begin{proof}
The cases involving \(1\) follow from the identity, and the case \(a=b\) from preservation of powers. By duality let \(ab=a\). Then \(a^i b^j=a^i\) for \(i,j>0\). Apply Lemma~\ref{lem:union} to \(C(a)\) on the left and \(C(b)\) on the right. We get
\begin{equation}\label{eq:absorption-intersection}
a\star b\in\bigl(C(a)\cup\{b\}\bigr)\cap\bigl(C(b)\cup\{a\}\bigr).
\end{equation}
If \(a^i=b^j\) with \(i,j>0\), then \(b^j=b^{j+1}\) since \(a^i b=a^i\). Thus \(e=b^j\) is idempotent, \(ae=a\) and \(eb=e\). Preservation of idempotent multiplication gives \(a\star b=a\star(e\star b)=a\star e=a\). Otherwise \eqref{eq:absorption-intersection} allows only \(a,b\), since \(a,b\) have no nontrivial product equal to \(1\) in \(D\).

Suppose \(a\star b=b\). We first exclude the case in which a power of \(b\) is a unit of a local monoid \(eMe\). If \(b^n\) is such a power, then \(ae=a\) since \(a b^n=a\), and \(b\) commutes with \(e\). For the latter, write \(v\) for the inverse of \(b^n\) in \(eMe\) and use \(be=b^n b v\) and \(eb=v b b^n\) to get \(be=ebe=eb\). Now \((eb)^n=b^n\), so \(eb\) is a unit of \(eMe\). It follows that \(a\star b=a\star(eb)=a(eb)=a\), by Proposition~\ref{prop:group-action} and preservation of idempotent multiplication. The dual reasoning in the \(\star\)-law excludes a positive power of \(a\) in a subgroup. Explicitly, suppose that \(a^n\), \(n\ge1\), is a unit of a local monoid \(eMe\); this condition does not depend on the law. In the \(\star\)-law \(a^n\star b=b\), hence \(e\star b=(e\star a^n)\star b=b\). As in the previous case, \(a\star e=e\star a\), and \(g:=a\star e=ae\) is a unit of \(eMe\). Proposition~\ref{prop:group-action} and \(eb=e\star b=b\) then give
\[
b=a\star b=(a\star e)\star b=g\star b=gb=a(eb)=ab=a,
\]
which is impossible. These properties of powers transport because powers, local monoids, and regular singletons have been identified. In particular both \(C(a)\) and \(C(b)\) are infinite in the remaining case, since a finite cyclic monoid has a group of eventual powers.

Put \(z=ba\). Then \(az=a^2\) and \(zb=z\). If \(z\) is a positive power of \(b\), the second equality contradicts distinctness of powers. If \(z\) is a positive power of \(a\), the first gives \(z=a\). In that case \(\{a\}\) belongs to the family of pointwise left fixed sets for \(b\). In the \(\star\)-law, a point fixed on the left by \(b\) is fixed on the left by \(a\), because \(a\star b=b\). The induced inclusion of fixed-set families transports back by Lemma~\ref{lem:fixed-sets}, forcing \(a^2=a\). This, too, is impossible.

We are left with \(z\notin C(a)\cup C(b)\). Put \(S=\langle1,a,b\rangle\), the same labelled submonoid in both laws by Lemma~\ref{lem:submonoids}. Since \(C(b)C(a)\) is transported pointwise, \(z=b^i\star a^j\) for some \(i,j>0\). In the \(\star\)-law the submonoid \(E=\{1\}\cup(S\star z)\) satisfies \(X\star E=E\cup X\) for \(X\in\Pow_1(S)\), because \(S\star z\) is a left ideal of \(S\). By Lemma~\ref{lem:union}, \(U_aE=E\cup U_a\) in the source. As \(a^2=az\in U_aE\) and \(a^2\notin U_a\), we have \(a^2\in E\setminus\{1\}\). Hence \(a^2=s\star z\) for some \(s\in S\). The relation \(a\star b=b\) puts each word in \(S\) in the form \(b^k\star a^\ell\), \(k,\ell\ge0\), by deleting \(a\) at occurrences of the adjacent pair \(ab\). Multiplication of \(a^2=s\star z\) on the right by \(b\) then gives
\[
b=(s\star z)\star b=s\star b^{i+1}=b^{k+i+1},
\]
contrary to distinctness. This proves \(a\star b=a\).
\end{proof}

\begin{lemma}\label{lem:pure}
If \(a,b\) lie in a common \(D_R\) or \(D_L\) and \(ab\in C(a)\cup C(b)\), then \(a\star b=ab\).
\end{lemma}

\begin{proof}
The case \(ab=1\) is trivial within \(D\). The cases \(ab=a,b\) follow from Lemma~\ref{lem:absorption}. By duality assume \(ab=a^r\) with \(r\ge2\) and \(ab\notin\{a,b\}\). For \(i\ge1\) and \(j\ge0\) we have
\begin{equation}\label{eq:pure-shift}
a^i b^j=a^{i+(r-1)j}.
\end{equation}
Consequently \(C(a)X=C(a)\cup X\) for \(X\in\Pow_1(C(b))\), and Lemma~\ref{lem:union} yields \(a\star b\in C(a)\cup\{b\}\). By Lemma~\ref{lem:absorption} in the reverse direction, \(a\star b\) is distinct from \(a,b\), so we can write \(a\star b=a^s\) with \(s\ge2\).

If \(C(a)\) is finite and both outputs belong to its periodic group, with identity \(e\), then \(e\) fixes both outputs on the left in their respective laws. Since \(ea\) belongs to a subgroup, preservation of subgroup multiplication gives
\[
a\star b=e\star(a\star b)=(ea)\star b=(ea)b=e(ab)=ab.
\]
Otherwise, interchanging the laws if necessary, the power \(a^r\) appears at a unique exponent. (If \(a^r\) lies in the periodic group of \(C(a)\), then \(a^s\) does not, and we exchange the roles of the two laws, that is, of \(\Psi\) and \(\Psi^{-1}\), together with those of \(r\) and \(s\).) Introduce the subsets
\[
Y_-=\{1,a\}\cup\{a^i:i>r\},\qquad Y_+=Y_-\cup\{a^r\}.
\]
We have \(Y_-U_b=Y_+U_b\), since \(ab=a^r\) and \(a^r b=a^{2r-1}\in Y_-\). By the membership test \eqref{eq:tail-membership}, \(\Psi(Y_+)\) contains \(a^r\), while \(\Psi(Y_-)\) omits it. Thus some point of \(\Psi(Y_-)\), multiplied by \(b\) in the other law, yields \(a^r\). This point is different from \(1\) because \(b\ne a^r\). Among the positive powers in \(\Psi(Y_-)\), only \(a\) occurs at an exponent at most \(r\), by the same membership test. Formula \eqref{eq:pure-shift} in the other law has \(s\) in place of \(r\). As \(a^r\) occurs at a unique exponent, the only way to obtain it is from \(a\star b=a^s\) with \(s=r\).
\end{proof}

\section{Finite cycles of multiplication}\label{sec:cycles}

We use the Green preorder \(x\leq_{\RR}y\) given by \(xM\subseteq yM\), or equivalently \(x=ys\) for some \(s\in M\). The left-handed preorder uses \(Mx\subseteq My\). Equality of ideals gives the relations \(\RR,\LL\), respectively. Left multiplication preserves \(\leq_{\RR}\), since \(x=ys\) implies \(tx=(ty)s\). Regularity is constant on either Green class, as in Proposition~\ref{prop:regular}.

\begin{lemma}\label{lem:finite-cycle}
If \(t^m x=x\) for some \(m\ge1\), then \(t\star x=tx\). Dually, \(xt^m=x\) implies \(x\star t=xt\).
\end{lemma}

\begin{proof}
We prove the left-handed statement. Put \(O=\{t^i x:i\ge0\}\) and \(A=C(t)\). Left multiplication by \(t\) is a cycle on \(O\). If \(O\) meets \(A\), the entire cycle lies in \(A\) and the result follows from preservation of powers. If \(x\) is regular, all of \(O\) is \(\LL\)-equivalent to \(x\), hence \(tx\) is regular and Proposition~\ref{prop:regular} applies. If \(t\) is a unit, Proposition~\ref{prop:group-action} gives the result. We can therefore assume \(O\) is nonregular, disjoint from \(A\), and \(t\) is a nonunit. Then \(t\) and \(O\) lie in one common \(D_R\) or \(D_L\), since all the points of \(O\) have neither one-sided inverse.

\medskip\noindent
\emph{Recognising the underlying cycle.}
If \(u,v\in O\) are comparable by \(\leq_{\RR}\), they are \(\RR\)-equivalent. Indeed if \(v=t^j u\leq_{\RR}u\), repeated application of \(t^j\) returns to \(u\) in the finite cycle, proving the reverse inequality. Let
\[
N=\bigcup_{k\ge2}O^k.
\]
We have \(N\cap O=\varnothing\). Otherwise \(w=o_1\cdots o_k\in O\) with \(k\ge2\). Write \(w=uv\) with \(v=o_k\). Since \(w\leq_{\RR}u\leq_{\RR}o_1\) and \(w\RR o_1\), we get \(w\RR u\). Also \(w\LL v\), since they lie on the same left cycle. Thus \(u=uvs\) and \(v=ruv\) for some \(r,s\in M\), giving \(u=uru\). This would make \(o_1\RR u\) regular. Also \(tN=N\), since multiplication by \(t\) on the first factor permutes \(O\).

For \(B\subseteq O\), \(E_B=\{1\}\cup N\cup B\) is a submonoid because any product of two nonidentity points lands in \(N\). For each \(y\in O\),
\begin{equation}\label{eq:cycle-set-products}
A E_{\{y\}}=A\cup N\cup O,\qquad A E_\varnothing=A\cup N.
\end{equation}
Using Lemma~\ref{lem:two-products} in the common \(D\), their products in the other law have these same point labels. Moreover \(A\cup N=A\star E_\varnothing\) is invariant under \(t\) on the left there. The difference between the two products in \eqref{eq:cycle-set-products} is exactly the set of points reached from \(y\) by left powers of \(t\) outside \(A\cup N\). Thus every \(y\in O\) reaches every other point of \(O\) in the target. If \(|O|=1\) the product is already preserved by Lemma~\ref{lem:absorption}. Otherwise mutual reachability gives a cycle, and its points are exactly \(O\) since an iterate entering \(A\cup N\) could never return to \(O\).

We can now apply the disjointness argument in the other law. As \(E_O\) is the submonoid generated by \(O\) in both laws and neither law has products of two or more points of \(O\) in \(O\), the subset \(N\) is the union of such products in both laws. In particular it is \(t\)-invariant on the left in both.

\medskip\noindent
\emph{Determining the direction of the cycle.}
Let \(\sigma,\sigma'\) be the left actions of \(t\) on \(O\) in the original and target laws, and put \(m=|O|\). If \(t\) has finitely many powers, take an idempotent power \(e\) of \(t\). Both actions of \(e\) are identities on \(O\), since an idempotent permutation is the identity. The action of \(t\) therefore agrees with that of \(te\) on \(O\) in both laws, and \(te\) lies in the cyclic group of eventual powers. Proposition~\ref{prop:group-action} gives the result. For \(m=2\) the cycles necessarily agree. Suppose \(m\ge3\) and \(C(t)\) is infinite. For \(r=1,2\) put
\[
F_r=\{t^i:i\ge0,\ i\not\equiv r\pmod m\}.
\]
These subsets are preserved pointwise by Lemma~\ref{lem:powers}. Let \(E=E_O\). For \(F=F_r\), the equations on \(X\in\Pow_1(D)\)
\begin{equation}\label{eq:cycle-families}
XE=FE,\qquad XE_\varnothing=X
\end{equation}
describe exactly the subsets \(X=F\cup N\cup Y\) for \(Y\subseteq O\). In fact \(FE=F\cup N\cup O\) because \(F\) contains \(1\) and each of its points preserves \(N,O\) on the left. Once \(X\subseteq F\cup N\cup O\), right multiplication by \(E\) adds \(N\cup O\) by union. Right multiplication by \(E_\varnothing\) adds \(N\) by union. The calculation holds in both laws.

For such \(X\) and \(z\in O\),
\[
XE_{\{z\}}=XE\quad\Longleftrightarrow\quad \sigma^r(z)\in Y.
\]
Indeed, the part added to \(Y\) by \(XE_{\{z\}}\) is \(F_r z=O\setminus\{\sigma^r(z)\}\), because \(F_r\) omits one residue class of exponents modulo \(m\). Other products with \(z\) contribute only to \(N\). The same test in the target gives the point \((\sigma')^r(z)\). Since \eqref{eq:cycle-families} transports, the part in \(O\) of the image of \(F_r\cup N\cup\{y\}\) is therefore
\[
\{(\sigma')^r\sigma^{-r}(y)\}.
\]
Finally \(U_tF_1=A=U_tF_2\), since every omitted exponent is one more than an included exponent. The products \(U_t(F_r\cup N\cup\{y\})\) are consequently equal for \(r=1,2\). Their image products have the same unordered pair \(\{z,\sigma'z\}\) in \(O\). Such a pair determines \(z\) for a cycle of length at least three, since otherwise it would be a two-cycle. This gives
\[
\sigma'\sigma^{-1}=(\sigma')^2\sigma^{-2},
\]
from which \(\sigma'=\sigma\). The dual assertion follows by reversing the products.
\end{proof}

\section{The mixed-product problem}\label{sec:mixed}

The preceding sections recover products whenever one factor belongs to a subgroup, whenever the result lies in a pure cyclic submonoid, or whenever a relevant translation has a finite cycle. The remaining case is the core of the proof: the subsets \(C(a^n)C(b^m)\) determine a lattice of congruences on the mixed-power orbit, but a congruence alone does not distinguish a translation from its inverse. We therefore reconstruct first the orbit and its origin, then its line distances, and finally the orientation of the two commuting translations.

Throughout this section \(a,b\) belong to a common \(D_R\) or \(D_L\). Put \(A=C(a)\), \(B=C(b)\), \(c=ab\) and \(d=a\star b\). Lemma~\ref{lem:pure}, applied in both directions, deals with the case in which either \(c\) or \(d\) belongs to \(A\cup B\). We therefore assume \(c,d\notin A\cup B\). In particular \(a,b\ne1\). Set
\[
O=\{a^i b^j:i,j\ge1\},\qquad
\alpha(x)=ax,\qquad \beta(x)=xb\qquad(x\in O).
\]
The maps \(\alpha,\beta\) commute by associativity. Every point of \(O\) is of the form \(\alpha^i\beta^j(c)\) with \(i,j\ge0\). We use \(O',\alpha',\beta'\) for the corresponding objects for the law \(\star\) and the initial point \(d\).

We repeatedly use the following observation. If \(r,s\ge0\), then
\begin{equation}\label{eq:recurrence}
c=\alpha^r\beta^s(c) \quad\Longrightarrow\quad
\alpha^r\beta^s=\id_O .
\end{equation}
Indeed, the equality at \(c\) extends to every point of \(O\) by commutativity. If \(r>0\), then \(\alpha^{r-1}\beta^s\) is an inverse of \(\alpha\). If \(s>0\), then \(\alpha^r\beta^{s-1}\) is an inverse of \(\beta\). Conversely, surjectivity of \(\alpha\) gives a preimage of \(c\) in \(O\) and hence a relation with \(r>0\). The corresponding assertion for \(\beta\) gives one with \(s>0\). For \(n\ge1\) write
\[
E_n(t)=\{1\}\cup\{t^i:i\ge n\}.
\]
These are submonoids, with the same point labels in the two laws.

\subsection{Finite orbits and noninvertible translations}

\begin{lemma}\label{lem:bounded-orbit}
With the notation above, \(c=d\) if \(O\) is finite or if at least one of \(\alpha,\beta\) is noninvertible.
\end{lemma}

\begin{proof}
Suppose first that neither map is invertible. By \eqref{eq:recurrence}, \(c\) is outside \(\alpha O\cup\beta O\), and every other mixed product belongs to this union. Consequently
\[
AB\setminus\bigl(A\cup B\cup E_2(a) B\cup A E_2(b)\bigr)=\{c\}.
\]
Indeed the two products involving \(E_2\) contain all mixed terms with at least one exponent greater than \(1\). Lemma~\ref{lem:two-products} identifies the displayed set pointwise in both laws. In the other law the only possible element of the difference is \(d\), so \(c=d\).

Next suppose that \(\alpha\) is invertible and \(\beta\) is noninvertible. A relation in \eqref{eq:recurrence} with \(r>0\) must have \(s=0\), so \(\alpha\) has finite order. Then
\[
K=O\setminus\beta O=\{\alpha^i(c):i\ge0\}.
\]
Every point outside the cycle on the right lies in \(\beta O\). If a point of the cycle belonged to \(\beta O\), a positive power of \(\alpha\) would take it back to \(c\), forcing \(c\in\beta O\) and hence making \(\beta\) invertible. Also \(K\cap A=\varnothing\), since \(A\) is closed under left multiplication by \(a\) and \(c\notin A\). If \(b\in K\), Lemma~\ref{lem:finite-cycle} gives \(a\star b=ab\). A point \(b^j\) with \(j\ge2\) cannot belong to \(K\). Otherwise some positive power of \(\alpha\) would give \(c=a^ib^j\in\beta O\). In the remaining case \(K\) is disjoint from \(A\cup B\), and
\begin{equation}\label{eq:layer}
K=AB\setminus\bigl(A\cup B\cup A E_2(b)\bigr).
\end{equation}
The right-hand side is determined in both laws. In the other law it is contained in \(O'\setminus\beta'O'\). Since it is nonempty, \(\beta'\) is noninvertible. Thus \(d\notin\beta'O'\) by \eqref{eq:recurrence} in the other law, and \(d\notin A\cup B\) by assumption. It follows that \(d\in K\).

Left multiplication by \(a\) has the same action on \(K\) in both laws by Lemma~\ref{lem:finite-cycle}. If \(m=|K|\), the product \(C(a^m)C(b)\) meets \(K\) in the singleton \(\{a^mb\}\). A point of the intersection is mixed since \(K\cap(A\cup B)=\varnothing\). Its \(b\)-exponent is \(1\) by \eqref{eq:layer}, and all positive \(a\)-exponents divisible by \(m\) act in the same way on the cycle. The same reasoning gives \(\{a^m\star b\}\) in the other law, since \(d\in K\) and the actions on \(K\) agree. The products of the cyclic submonoids have the same point labels. One further common left multiplication by \(a\) takes these two points to \(c,d\), respectively. If \(\beta\) is invertible and \(\alpha\) is noninvertible, apply this argument in the opposite monoids with \(a,b\) interchanged.

It remains to consider the case in which \(O\) is finite and both maps are invertible. By Lemma~\ref{lem:finite-cycle}, both actions on \(O\) are preserved. Since \(AB\) is determined pointwise and \(d\notin A\cup B\), we have \(d\in O\). The permutations \(\alpha,\beta\) act transitively together on \(O\), and each inverse is a positive power. Thus their nonnegative powers applied to \(d\) also give all of \(O\), proving \(O'=O\).

Suppose \(O\) meets \(A\) at \(a^k\). A positive power of \(\alpha\) fixes \(a^k\), giving \(a^{k+i}=a^k\) for some \(i>0\). Thus \(a\) has finitely many powers. Let \(e\) be the identity of the group of its eventual powers. It is a positive power of \(a\), so its left action is an idempotent permutation of \(O\), hence the identity. This holds in both laws. The element \(ea\) belongs to that subgroup, so Proposition~\ref{prop:group-action} gives
\[
c=(ea)b=(ea)\star b=d.
\]
If \(O\) meets \(B\), use the same reasoning with right multiplication. If \(O\) is disjoint from \(A\cup B\), let \(m,n\) be the orders of \(\alpha,\beta\). The intersection \(C(a^m)C(b^n)\cap O\) is the singleton \(\{a^mb^n\}\). All its mixed terms have exponents divisible by \(m,n\), which makes them equal. The corresponding singleton in the other law is the same by Lemma~\ref{lem:two-products}. Applying the two common translations once each sends that point to \(c\) in the original law and to \(d\) in the other law, because \(\alpha^m=\beta^n=\id_O\).
\end{proof}

\subsection{The group of translations of an infinite orbit}

The point of passing to a translation group is that a mixed word becomes a lattice point, while products of cyclic submonoids become cosets of explicitly described subgroups. This converts the reconstruction of \(ab\) into a question about the geometry of an abelian group.

In a possible disagreement we may now assume, by Lemma~\ref{lem:bounded-orbit} in both directions, that \(O,O'\) are infinite and both translations are invertible in each law. In this situation \(\alpha,\beta\) have infinite order. Indeed, a relation in \eqref{eq:recurrence} with \(s>0\) would give \(\beta\) finite order if \(\alpha\) had finite order. The orbit would then be finite. Interchanging the roles of the maps proves the other assertion. In particular \(a\) and \(b\) each have distinct powers. Adding a relation with \(r>0\) to one with \(s>0\) gives a relation with both coefficients positive.

Let \(G\) be the abelian group generated by the commuting permutations \(\alpha,\beta\) of \(O\), and write \(u,v\) for the respective generators in additive notation. The group acts transitively on \(O\). A permutation in \(G\) which fixes one point fixes all points by commutativity and transitivity. Hence the action is free, that is, for any \(x,y\in O\) exactly one element of \(G\) sends \(y\) to \(x\). We denote that element by \(x-y\). Choose \(o=(\alpha\beta)^{-1}(c)\) and write \(o+g\) for the action of \(g\in G\) on \(o\). Then
\begin{equation}\label{eq:grid-coordinates}
a^i b^j=o+iu+jv\qquad(i,j\ge1).
\end{equation}

Let \(q>0\) be least with \(qv\in\ZZ u\). There is a unique \(p>0\) with \(pu+qv=0\). To see the sign and uniqueness, use the relation with both coefficients positive and the fact that \(u\) has infinite order. Dividing the coefficient of \(v\) in any integer relation by \(q\) shows that every relation is a multiple of \((p,q)\). Explicitly, if \(iu+jv=0\), write \(j=kq+j'\) with \(0\le j'<q\); then \(j'v\in\ZZ u\) forces \(j'=0\) by the minimality of \(q\), and \((i-kp)u=0\) forces \(i=kp\). For the relation with both coefficients positive this gives \(k>0\), hence \(p>0\). In particular
\begin{equation}\label{eq:translation-group}
G\cong\ZZ^2/\ZZ(p,q),\qquad
[G:\ZZ u]=q,\qquad [G:\ZZ v]=p .
\end{equation}
With \(h=\gcd(p,q)\), the homomorphism \(G\to\ZZ\) given by \([(i,j)]\mapsto qi-pj\) has kernel generated by \([(p/h,q/h)]\), of order \(h\). It follows that
\begin{equation}\label{eq:residually-divisible-intersection}
\bigcap_{n\ge1} nG=\{0\}.
\end{equation}
Indeed any element of the intersection has zero image in \(\ZZ\). It is also \(hy\) for some \(y\in G\), and \(y\) must have zero image in \(\ZZ\). Since the kernel has order \(h\), we get \(hy=0\).

For \(n,m\ge1\) set \(H_{n,m}=\langle nu,mv\rangle\le G\). Every element of this subgroup is a combination of \(nu,mv\) with strictly positive integer coefficients. Starting with integer coefficients, add a sufficiently large multiple of the relation
\begin{equation}\label{eq:positive-multiple}
(mp)(nu)+(nq)(mv)=0.
\end{equation}
This positivity will allow the products of submonoids to record the full subgroups \(H_{n,m}\).

\subsection{Recovering the orbit and its origin}

\begin{lemma}\label{lem:orbit-origin}
The mixed orbits \(O,O'\) have the same point labels, and their origins \(o,o'\) coincide.
\end{lemma}

\begin{proof}
If \(a^k\in O\) for \(k>0\), then applying \(\beta\) gives \(a^k=o+ku\) by \eqref{eq:grid-coordinates}. Applying \(\alpha\) proves \(b^\ell=o+\ell v\) when \(b^\ell\in O\) and \(\ell>0\). In particular \(A\cap B=\{1\}\). Indeed, \(a^k=b^\ell\) with \(k,\ell>0\) implies \(a^{k+1}=ab^\ell\in O\), hence \(ku=\ell v\). This is impossible since all nonzero relations in \eqref{eq:translation-group} have coefficients of the same sign.

For all \(n,m\ge1\), the mixed terms of \(C(a^n)C(b^m)\) cover \(o+H_{n,m}\) by \eqref{eq:grid-coordinates} and \eqref{eq:positive-multiple}. Any pure positive power from that product lying in \(O\) already belongs to that coset by the coordinate formulas for pure powers. We have \(1\notin O\) since \(a,b\) are nonidentity elements in a common \(D_R\) or \(D_L\). We conclude that
\begin{equation}\label{eq:orbit-products}
C(a^n)C(b^m)\cap O=o+H_{n,m}.
\end{equation}
In particular for \(k>0\),
\begin{equation}\label{eq:pure-in-orbit}
a^k\in O \quad\Longleftrightarrow\quad
a^k\in C(a^n)C(b)\ \text{for some }n>k .
\end{equation}
For the forward implication take \(n>k\) coprime to \(p\), since then \(nu\) generates \(G/\ZZ v\) and \(H_{n,1}=G\). For the converse, \(a^k\) belongs to neither \(C(a^n)\) nor \(C(b)\), since the powers of \(a\) are distinct and \(A\cap B=\{1\}\). Any such representation is therefore mixed. Dually, \(b^\ell\) with \(\ell>0\) belongs to \(O\) if and only if it belongs to \(C(a)C(b^m)\) for some \(m>\ell\).

These membership tests transport by Lemma~\ref{lem:two-products}. Also \(O\setminus(A\cup B)=AB\setminus(A\cup B)\), which is determined pointwise. This proves \(O'=O\). Since \(H_{n,n}=nG\), equations \eqref{eq:orbit-products} and \eqref{eq:residually-divisible-intersection} now show that
\[
\bigcap_{n\ge1}\bigl(C(a^n)C(b^n)\cap O\bigr)=\{o\}.
\]
This is the same intersection in both laws.
\end{proof}

\subsection{Recovering congruences and distances}

We now show that multiplication in \(\Pow(M)\) detects which orbit points are congruent modulo \(H_{n,m}\). Intersecting these congruences leaves the two families of cyclic lines; varying the modulus then recovers the integer distance on each line. Importantly, this step identifies distances but not yet the positive direction.

\begin{lemma}\label{lem:line-distances}
The families of subsets \(x+\ZZ u\) and \(x+\ZZ v\) of \(O\) are determined in both laws. Their numbers \(q,p\), respectively, and the distances along each such subset are also determined. The distance from \(x\) to \(x+ku\) on a line in the \(u\)-direction is \(|k|\), with the analogous convention in the \(v\)-direction.
\end{lemma}

\begin{proof}
For \(x,y\in O\) we have the test
\begin{equation}\label{eq:coset-test}
C(a^n)U_xC(b^m)=C(a^n)U_yC(b^m)
\quad\Longleftrightarrow\quad x-y\in H_{n,m}.
\end{equation}
The two sides have the common part \(C(a^n)C(b^m)\) contributed by the identity in \(U_x,U_y\). Its intersection with \(O\) is \(o+H_{n,m}\) by \eqref{eq:orbit-products}. The extra parts are
\[
C(a^n)xC(b^m)=x+H_{n,m},\qquad
C(a^n)yC(b^m)=y+H_{n,m},
\]
because both translations preserve \(O\) and \eqref{eq:positive-multiple} yields all exponents needed to cover a coset. Outside \(O\) the two sides are the same. Inside \(O\) equality holds precisely when the two added cosets coincide. Thus \eqref{eq:coset-test} transports through \(\Psi\).

At fixed \(n\), the intersection of all \(H_{n,m}\) over \(m\ge1\) is \(n\ZZ u\). Every term contains \(n\ZZ u\), and for \(m=nq\) we have \(mv=-npu\), giving equality. At fixed \(m\), the intersection over \(n\ge1\) is \(m\ZZ v\), using \(n=mp\). The tests defining these intersections run over every positive \(m\) or \(n\). Taking \(n=1\) or \(m=1\) therefore identifies the two families of lines before \(q,p\) are counted. The numbers of their classes give \(q,p\) by \eqref{eq:translation-group}.

If \(y-x=ku\) with \(k\ne0\), then \(y-x\in n\ZZ u\) holds precisely for positive divisors \(n\) of \(k\). Their largest element is \(|k|\). This determines distances on the \(u\)-lines, and the same argument applies to the \(v\)-lines. As a result, on every line, its two integer coordinate systems in the two laws differ by an isometry of \(\ZZ\) with the usual distance.
\end{proof}

\subsection{Recovering the directions of the translations}

The remaining ambiguity is a reflection on each cyclic line. An intrinsic test involving a tail of powers distinguishes the positive half-line from the negative one. Two consecutive tails then identify the first positive step, which is enough to remove the reflection ambiguity.

\begin{lemma}\label{lem:line-directions}
Under the infinite-orbit assumptions, \(c=d\).
\end{lemma}

\begin{proof}
\emph{Step 1: an intrinsic positive ray.} For \(n\ge1\), define
\[
T_n=\{x\in O:E_n(a)U_x U_{b^q}=E_n(a)U_{b^q}\}.
\]
The set \(T_n\) is intrinsic: the orbit \(O\) and the integer \(q\) have been recovered in Lemmas~\ref{lem:orbit-origin} and~\ref{lem:line-distances}, and \(E_n(a)\) and \(U_{b^q}\) are preserved by Lemma~\ref{lem:powers}. We claim that
\begin{equation}\label{eq:orientation-test}
T_n\setminus\{o\}=\{o+ku:k\ge n\}.
\end{equation}
To establish the claim, put \(S_n=E_n(a)U_{b^q}\). Its mixed terms are
\[
a^ib^q=o+iu+qv=o+(i-p)u\qquad(i\ge n),
\]
where \(qv=-pu\). Thus they fill the ray \(\{o+ju:j\ge n-p\}\). Any pure power of \(a\) in \(S_n\cap O\) is of the form \(o+iu\) with \(i\ge n\), and hence already lies on this ray. Apart from the ray, \(S_n\cap O\) may contain only \(b^q=o-pu\). Notice that \(o-pu\) is strictly outside the ray, since \(n\ge1\).

If \(x\in T_n\), then, because \(1\in E_n(a)\), both \(x\) and \(xb^q\) must lie in \(S_n\cap O\). In particular \(x=o+ku\) for some \(k\in\ZZ\), and \(xb^q=o+(k-p)u\). If \(xb^q\) lies on the ray, then \(k-p\ge n-p\), so \(k\ge n\). The only other possibility is \(xb^q=b^q\), which forces \(x=o\) by the free translation action. This proves one inclusion in~\eqref{eq:orientation-test}.

Conversely, if \(x=o+ku\) with \(k\ge n\), the terms newly introduced in \(E_n(a)U_xU_{b^q}\) are
\[
x,\qquad xb^q,\qquad a^ix,\qquad a^ixb^q\qquad(i\ge n).
\]
Their respective \(u\)-coordinates relative to \(o\) are \(k\), \(k-p\), \(k+i\) and \(k+i-p\). Each is at least \(n-p\); all these terms therefore belong to \(S_n\). The required equality follows and the claim is proved.

\emph{Step 2: identifying the first positive steps.} Equation~\eqref{eq:orientation-test} gives an intrinsic singleton,
\[
T_1\setminus(T_2\cup\{o\})=\{o+u\}.
\]
The symmetric test
\[
\{x\in O:U_{a^p}U_xE_n(b)=U_{a^p}E_n(b)\}
\]
selects \(\{o+kv:k\ge n\}\), up to the possible exceptional point \(o\). It therefore identifies \(o+v\) in the same fashion. Consequently both \(o+u\) and \(o+v\) are the same labelled points in the two laws. Since distances on the \(u\)- and \(v\)-lines have already been recovered, these two points determine the positive directions on the respective lines through \(o\).

\emph{Step 3: propagation to all parallel lines.} Every \(u\)-line meets the \(v\)-line through \(o\), because \(u\) and \(v\) generate \(G\). Let \(y\) be an intersection point. The \(v\)-line has the same metric and the same positive direction in the two laws, so the point
\[
y-qv=y+pu
\]
is intrinsically identified. The two distinct points \(y\) and \(y+pu\) determine the positive direction of the \(u\)-line through \(y\): an isometry of \(\ZZ\) fixing two distinct points is the identity. Hence \(\alpha\) and \(\alpha'\) agree on every \(u\)-line. Likewise \(\beta(o)=\beta'(o)\) by Step~2. We conclude that
\[
c=\alpha\beta(o)=\alpha'\beta'(o)=d.
\]
\end{proof}

\begin{proposition}\label{prop:common-D}
For any \(a,b\) in a common \(D_R\) or \(D_L\), \(a\star b=ab\).
\end{proposition}

\begin{proof}
Lemma~\ref{lem:pure} covers a product in \(A\cup B\) in either law. Lemma~\ref{lem:bounded-orbit}, applied in either direction, covers a finite mixed orbit or a noninvertible translation. Lemma~\ref{lem:line-directions} covers the remaining case.
\end{proof}

\subsection{An orbit with torsion}

The parameters \(p\) and \(q\) in~\eqref{eq:translation-group} need not be coprime. The following example illustrates why it is necessary to use the whole family of coset tests rather than only the rank-one quotient of the translation group.

\begin{example}\label{ex:torsion}
Consider the monoid with presentation
\[
M=\langle a,b\mid a^3b^5=ab\rangle_{\mathrm{mon}}.
\]
Neither side of the defining relation is empty, so a nonempty word cannot represent \(1\). In particular \(a,b\) have neither a left nor a right inverse. Inside a word \(a^i b^j\), \(i,j\ge1\), replacing either side of the relation by the other preserves the form of the word and changes the exponent pair by \(\pm(2,4)\). Conversely, any two positive exponent pairs differing by a multiple of \((2,4)\) are joined by these replacements. Start from the smaller pair and use the replacement \(ab\mapsto a^3b^5\) repeatedly. Pure powers cannot occur in such a chain, so they are disjoint from \(O\). Every class modulo \(\ZZ(2,4)\) has a positive representative, since adding a large multiple of \((2,4)\) makes both coordinates positive. This gives
\[
O\simeq\ZZ^2/\ZZ(2,4),\qquad
u=[(1,0)],\quad v=[(0,1)],\quad o=[(0,0)].
\]
Here the origin has the representative \(a^2b^4\). The coordinate function \([(i,j)]\mapsto 2i-j\) (that is, \((qi-pj)/2\), since here \(p=2\) and \(q=4\)) takes the same value on \(o\) and \(o+u+2v\), which are distinct because \(u+2v\) has order two. The set \(C(a^2)C(b^4)\cap O\) distinguishes the two points. It equals \(o+\ZZ(2u)\), since \(4v=-2u\), and contains only \(o\) of the two points under consideration. This illustrates how the coset tests recover the torsion coordinate.
\end{example}

\section{Completion of the reconstruction}\label{sec:completion}

We combine the preceding reconstruction results. The only pairs not already covered are the mixed one-sided-inverse configurations; these reduce to products for which regularity or the common-\(D\) argument is available.

\begin{proof}[Proof of Theorem~\ref{thm:main}]
Section~\ref{sec:labels} gives the identity of \(T\) and the unique bijection \(f\) satisfying \eqref{eq:unit-label} and \eqref{eq:pair-label}. Lemma~\ref{lem:submonoids} and Proposition~\ref{prop:regular} give the assertions about submonoids and regular singletons. It remains to prove that \(ab=a\star b\) for every \(a,b\in M\).

Proposition~\ref{prop:group-action} covers products with a unit. Proposition~\ref{prop:common-D} covers pairs in a common \(D_R\) or \(D_L\). Thus any remaining pair consists of one nonunit with a left inverse and one nonunit with a right inverse. Indeed, if \(a,b\) are nonunits not lying in a common \(D_R\) or \(D_L\), then one of them has a right inverse and one of them has a left inverse; since a nonunit cannot have both, the two factors have one-sided inverses on opposite sides. In particular both factors are regular.

First suppose \(sa=1=bt\) for some \(s,t\). The element \(ab\) is regular, since
\[
(ab)(ts)(ab)=a(bt)(sa)b=ab .
\]
Its product is therefore preserved by Proposition~\ref{prop:regular}.

For the other order, suppose \(as=1=tb\) and put \(c=ab\). If \(c\) is regular the same proposition applies. Suppose \(c\) is nonregular. Then \(s,c\in D_R\). Indeed \(s\) has a left inverse and is a nonunit, and \(c\) has neither one-sided inverse. Moreover \(sc\LL c\) because \(sc\) is a left multiple of \(c\) and \(a(sc)=c\). Hence \(sc\) is nonregular as well, giving \(sc,t\in D_L\). We can now use Proposition~\ref{prop:common-D} twice to obtain
\[
s\star c\star t=(sc)t=(sa)(bt)=(sa)\star(bt).
\]
The last equality follows from preservation of idempotent multiplication, since \(sa,bt\) are idempotent. The products of regular singletons also give
\[
a\star s=1,\qquad t\star b=1,\qquad s\star a=sa,\qquad b\star t=bt.
\]
Multiplying the preceding product identity on the left by \(a\) and on the right by \(b\) in the \(\star\)-law now gives \(c=a\star b\). This proves multiplicativity and completes the proof.
\end{proof}

\section{Singletons and automorphisms}\label{sec:applications}

\subsection{Monoids generated by regular elements}

Theorem~\ref{thm:main} specifies how to recover all point labels. In some monoids those labels also determine the images of all singletons.

\begin{corollary}\label{cor:regular-generated}
Suppose \(M\) is generated as a monoid by its regular elements, and let \(\Phi:\Pow(M)\to\Pow(T)\) be an isomorphism for a semigroup \(T\). Then \(\Phi(\{a\})=\{f(a)\}\) for every \(a\in M\), where \(f\) is the isomorphism of Theorem~\ref{thm:main}.
\end{corollary}

\begin{proof}
Write \(a=r_1\cdots r_k\) with all \(r_i\) regular. The theorem gives
\[
\Phi(\{a\})
=\Phi(\{r_1\})\cdots\Phi(\{r_k\})
=\{f(r_1)\cdots f(r_k)\}
=\{f(a)\}.
\]
For the empty product this is preservation of the identity.
\end{proof}

In particular this applies to every regular monoid, hence to every inverse monoid. For inverse semigroups Liu and Tringali obtained singleton preservation even in the absence of an ambient identity~\cite{LiuTringali}.

\subsection{The induced map on automorphism groups}

For each isomorphism \(\Phi\) as in Theorem~\ref{thm:main}, write \(f_\Phi\) for the isomorphism it induces. Its two labelling identities imply
\[
f_{\Theta\circ\Phi}=f_\Theta\circ f_\Phi
\]
for any composable isomorphisms of power monoids. Indeed both sides agree on units by applying the first identity to the singleton of a unit, and on nonunits by applying the second identity to a nonunit pair. Conversely any monoid isomorphism \(h:M\to T\) gives an isomorphism of power monoids \(X\mapsto h[X]\), whose induced map is \(h\).

\begin{corollary}\label{cor:aut}
For every monoid \(M\), the map \(\Phi\mapsto f_\Phi\) is a split surjective homomorphism
\[
\Aut(\Pow(M))\longrightarrow\Aut(M).
\]
If \(K(M)\) is its kernel, then
\[
\Aut(\Pow(M))\cong K(M)\rtimes\Aut(M),
\]
with the action on \(K(M)\) given by conjugation by the setwise action of \(\Aut(M)\).
\end{corollary}

\begin{proof}
Compatibility with composition gives the homomorphism. The map \(h\mapsto(X\mapsto h[X])\) is its homomorphic section. The semidirect product decomposition follows from this section.
\end{proof}

\subsection{A family with an explicit kernel}

Let \(A\) be any set, and let \(M_A=\{1,0\}\sqcup A\) be the monoid with identity \(1\) and with \(xy=0\) whenever \(x,y\ne1\). Write \(\mathcal B(A)=\{B:\varnothing\ne B\subseteq A\}\), and \(\Sym(Y)\) for the group of all bijections of a set \(Y\).

\begin{proposition}\label{prop:null-kernel}
There are isomorphisms
\[
K(M_A)\cong\Sym(\mathcal B(A)),\qquad
\Aut(\Pow(M_A))\cong\Sym(\mathcal B(A))\rtimes\Sym(A),
\]
where \(\Sym(A)\) acts by conjugation via \(B\mapsto h[B]\) for \(h\in\Sym(A)\). In particular, if \(|A|=m<\infty\), then
\[
|\Aut(\Pow(M_A))|=(2^m-1)!\,m! .
\]
\end{proposition}

\begin{proof}
An automorphism of \(M_A\) fixes \(1\) and \(0\), since \(0\) is the only idempotent distinct from \(1\). Every permutation of \(A\) extends to such an automorphism. Thus \(\Aut(M_A)=\Sym(A)\) under restriction.

Let \(\Phi\in K(M_A)\). Theorem~\ref{thm:main} implies that all submonoids of \(M_A\), the pair \(U_0=\{1,0\}\), and the singleton \(\{0\}\) are fixed. We have \(1\in X\) if and only if \(XM_A=M_A\), so \(\Phi\) preserves the family of subsets containing \(1\). Every such subset which also contains \(0\) is a submonoid and is fixed. Within that family, the idempotent subsets are \(\{1\}\) and the subsets containing \(0\). It follows that a subset with \(1\in X\) and \(0\notin X\) is mapped to another with the same memberships in \(1,0\). Now for every nonempty \(X\),
\[
XU_0=U_0X=X\cup\{0\}.
\]
If \(1\in X\) and \(0\notin X\), the right side is fixed, so \(X\) itself is fixed.

Among the subsets omitting \(1\), the equation \(XU_0=X\) recognises those containing \(0\). Hence \(\Phi\) restricts to a permutation \(\pi\) of \(\mathcal B(A)\), and the same product with \(U_0\) forces
\[
\Phi(B)=\pi(B),\qquad \Phi(\{0\}\cup B)=\{0\}\cup\pi(B)
\qquad(B\in\mathcal B(A)).
\]
These formulas and the already fixed subsets determine \(\Phi\).

Conversely, given any permutation \(\pi\) of \(\mathcal B(A)\), the formulas define a permutation of the full power monoid if we fix \(\{0\}\) and every subset containing \(1\). To check products, if both factors contain \(1\), the product also contains \(1\), and the factors and the product are fixed. If both factors omit \(1\), their product is \(\{0\}\). If \(1\notin X\) and \(1\in Y\), then \(XY=YX=X\) for \(Y=\{1\}\) and \(XY=YX=X\cup\{0\}\) otherwise. The formulas commute with adjoining \(0\), which checks these products. This gives an automorphism in \(K(M_A)\). The correspondence with \(\pi\) respects composition, and conjugation by the section in Corollary~\ref{cor:aut} gives the stated action.
\end{proof}

For example, with \(A=\{a,b\}\), take the permutation of \(\mathcal B(A)\) exchanging \(\{a\}\) and \(\{a,b\}\) and fixing \(\{b\}\). Its extension belongs to \(K(M_A)\) and moves a singleton. Thus the isomorphism induced by a power-monoid isomorphism need not describe its action on nonregular singletons.

\begin{remark}\label{rem:kernel}
The kernel \(K(M)\) can be nontrivial also for groups: for \(G=\ZZ/3\) one has \(K(G)\cong\ZZ/3\) (Example~\ref{ex:units}). A computer search over all monoids of order at most \(5\) shows that \(K(M)\) is trivial for \(30\) of the \(35\) monoids of order \(4\) and for \(201\) of the \(228\) monoids of order \(5\). It would be interesting to characterise the monoids with trivial \(K(M)\).
\end{remark}

\section*{Declaration of AI-assisted research and manuscript preparation}

OpenAI's ChatGPT and Codex were used during the development and preparation of this manuscript. This assistance included exploration of the reconstruction strategy, formulation and preliminary checking of intermediate arguments, development of the mixed-orbit and coset tests, and revision of the exposition and \LaTeX{} source. The named authors are responsible for the content of the manuscript; the use of these tools does not constitute authorship.

\end{document}